\pdfoutput=1
\documentclass[11pt,reqno]{amsart}
\usepackage[paperheight=279mm,paperwidth=18cm,textheight=26cm,textwidth=14cm,includehead]{geometry}
\usepackage[T1]{fontenc}
\usepackage[utf8]{inputenc}
\usepackage[unicode]{hyperref}
\hypersetup{
bookmarks=true,
colorlinks=true,
citecolor=[rgb]{0,0,0.5},
linkcolor=[rgb]{0,0,0.5},
urlcolor=[rgb]{0,0,0.75},
pdfpagemode=UseNone,
pdfstartview=FitH,
pdfdisplaydoctitle=true,
pdftitle={Bi-parameter Carleson embeddings with product weights},
pdfauthor={Arcozzi, Mozolyako, Psaromiligkos, Volberg, Zorin-Kranich},
pdflang=en-US
}

\usepackage[style=alphabetic,maxalphanames=4]{biblatex}
\usepackage{amssymb}
\usepackage{amsmath}
\usepackage{amsthm}
\usepackage{amsfonts}
\usepackage{mathtools}

\usepackage{tikz}
\usetikzlibrary{patterns}

\DeclarePairedDelimiterX\Set[1]\{\}{%

#1
}

\DeclarePairedDelimiterX\innerp[2]{\langle}{\rangle}{#1,#2}

\numberwithin{equation}{section}
\theoremstyle{plain}
\newtheorem{theorem}[equation]{Theorem}

\newtheorem{lemma}[equation]{Lemma}

\newtheorem{conj}[equation]{Conjecture}

\theoremstyle{remark}
\newtheorem{remark}[equation]{Remark}

\theoremstyle{definition}

\newcommand{\one}{\mathbf{1}}

\DeclareMathOperator{\supp}{supp}

\newcommand{\pd}{\partial}

\newcommand{\al}{\alpha}
\newcommand{\eps}{\varepsilon}
\newcommand{\la}{\lambda}
\newcommand{\vf}{\varphi}
\newcommand{\om}{\omega}

\newcommand{\cE}{\mathcal E}
\newcommand{\cF}{\mathcal F}

\newcommand{\cT}{\mathcal T}

\newcommand{\bD}{\mathbb D}
\newcommand{\bE}{\mathbb E}

\newcommand{\bI}{\mathbb I}

\newcommand{\bR}{\mathbb R}

\newcommand{\bV}{\mathbb V}

\begin{document}

\title[The lack of small energy majorization in dimension $2$]%
{Differences between the potential theories on a tree and on a bi-tree}
\author[P.~Mozolyako]{Pavel Mozolyako}
\thanks{PM is supported by the Russian Science Foundation grant 17-11-01064}
\address[P.~Mozolyako]{St. Peterburg State University}
\email{pmzlcroak@gmail.com}
\author[A.~Volberg]{Alexander Volberg}
\thanks{AV is partially supported by the NSF grant DMS-160065 and DMS 1900268 and by Alexander von Humboldt foundation}
\address[A.~Volberg]{Department of Mathematics, Michigan Sate University, East Lansing, MI. 48823 and Hausdorff Center of Universit\"at Bonn}
\email{volberg@math.msu.edu}
%
%
\begin{abstract}
In this note we give several counterexamples. One shows that small energy majorization on bi-tree fails. The second counterexample shows that energy estimate $\int_T \bV^\nu_\eps \, d\nu \le C \eps |\nu|$ always valid on a usual tree by a trivial reason (and with constant $C=1$) cannot be valid in general on bi-tree with any $C$ whatsoever. On the other hand, a weaker estimate $\int_{T^2} \bV^\nu_\eps \, d\nu \le C_\tau \eps^{1-\tau} \cE[\nu]^{\tau} |\nu|^{1-\tau}$  is valid on bi-tree with any $\tau>0$. It is proved in \cite{MPVZ} and is called improved surrogate maximum principle for potentials on bi-tree. The estimate $\int_{T^3} \bV^\nu_\eps \, d\nu \le C_\tau \eps^{1-\tau} \cE[\nu]^{\tau} |\nu|^{1-\tau}$  with $\tau=2/3$ holds on tri-tree. We do not know {\it any} such estimate with {\it} any $\tau<1$ on four-tree. The third counterexample disproves the estimate $\int_{T^2} \bV^\nu_x \, d\nu \le F(x)$ for any $F$ whatsoever for some probabilistic $\nu$ on bi-tree $T^2$. On a simple tree $F(x)=x$ would suffice to make this inequality to hold.
\end{abstract}
\maketitle

\section{Introduction. Potential theory on multi-trees}
\label{intro}
Embedding theorems on graphs are interesting in particular because they are related to the structure of spaces of holomorphic functions. For Dirichlet space on  disc $\bD:=\{z: |z|<1\}$  this fact has been explored in  \cite{ARS2002}, \cite{ARSW},
\cite{ARSW11}, and for Dirichlet space on bi-disc $\bD^2$ in \cite{AMPS}, \cite{AMPVZ-K}, \cite{AHMV}. Bi-disc case is much harder as the corresponding graph has cycles. One particular interesting case is studied  in \cite{Saw1}, where a small piece of bi-tree is considered.

The difference between one parameter theory (graph is a tree) and two parameter theory (graph is a bi-tree) is huge.
One explanation is that in a multi-parameter theory all the notions of singular integrals, para-products, BMO, Hardy classes etc become much more subtle than in one parameter  settings. There are many examples of this effect.  It was demonstrated in results of S.Y. A. Chang, R. Fefferman and L. Carleson, see \cite{Carleson}, \cite{Chang}, \cite{ChF1},\cite{RF1}, \cite{TaoCar}.

The papers dealing with poly-disc and multi-trees mentioned above are all have a common feature: they are based on potential theory on multi-trees. Let us recall the reader the main notations and facts of such a theory. We will do this for bi-tree just for the sake of simplicity. 

Let $T$ denote the dyadic rooted tree with root $o$, we can associate the vertices with dyadic sub-intervals of $I^0:=[0,1]$, and $o$ with $I^0$ itself.  Similarly, let $T^2$ denote the dyadic rooted bi-tree with root $o$, we can associate the vertices with dyadic sub-rectangles of $Q^0:=[0,1]\times [0,1]$, and $o$ with $Q^0$ itself.  Both objects have partial order, which is the same as inclusion for intervals, rectangles correspondingly. 

Both objects have a natural integration operator, if $f$ is a non-negative function on $T$ or $T^2$, and $\al$ is a vertex of $T$ or $T^2$, then
$$
\bI f(\al) := \sum_{\al'\ge \al} f(\al')\,.
$$
We can call $\bI$ the Hardy operator on a corresponding graph: it sums up values from $\al$ to $o$ along all directed paths from $o$ to $\al$. For $T$ such a path is unique for any $\al$, for $T^2$ there are many such paths.

The formally adjoint operator is $\bI^*$ and
$$
\bI^* f(\al)  := \sum_{\al'\le \al} f(\al')\,.
$$
Let us make a convention that always our $T$ and/or $T^2$ are {\it finite graphs}, maybe very deep, but finite, and leaves are dyadic intervals  of size $2^{-N}$ in the case of $T$ or dyadic squares of size $2^{-N}\times 2^{-N}$ in the case of $T^2$.
Then $\bI^*$ is always defined. The set of leaves is a ``boundary'' of the graph and is denoted $\pd T$ or $\pd T^2$ correspondingly.

Now we want to introduce potential of measure. Again for simplicity (this is not at all important) let us call {\it measure} the function $\mu$ on $T^2$ that is identically zero on $T^2\setminus \pd T^2$ and just an arbitrary non-negative function on $\pd T^2$. We have the same way to define measure on $T$. Of course, what we really doing is defining {\it granular measures} on $Q^0$ and $I^0$ correspondingly. Here {\it granular} means that our measure have constant density  with respect to dyadic squares of size $2^{-N}\times 2^{-N}$ or dyadic intervals of size $2^{-N}$ correspondingly.

We wish to have all estimates ever met in our theory to {\it not depend} on $N$. Then by making limit when $N\to\infty$ we can consider {\it all} measures on $Q^0$ or $I^0$ eventually.

Given such a measure $\mu$ we define its potential at a vertex $\al$ of $T$ or $T^2$ as
$$
\bV^\mu (\al) := \bI\circ \bI^*(\mu) (\al)\,.
$$

Notice that as $\al$ is actually a dyadic rectangle $R=I\times J$ inside $Q^0$ (or dyadic interval $I$ inside $I^0$), then $\bI^*(\mu)(\al)$ is just $\mu(R)$ ($\mu(I)$ correspondingly).

\bigskip

But $\bV^\mu(\al)$ is a more complicated object, it is the sum of $\mu(R')$ over all $R'$ containing $R$, where $R$ is associated with vertex $\al\in T^2$ (correspondingly the sum of $\mu(I')$ over all $I'$ containing $I$, where $I$ is associated with vertex $\al\in T$).

\bigskip

Let us be on $T$ and let $\bV^\mu \le 1$ on $\supp\mu$ (these are vertices  of $\pd T$ where $mu>0$. Then we can easily see that $\bV^\mu\le 1$ everywhere. In fact, without loss of generality $\mu\neq 0$, and let $\beta\in \pd T$ and let $\mu(\beta)=0$. 

\medskip

We can find unique smallest predecessor $\gamma>\beta$ such that there is $\al\in \pd T$, $\mu(\al)>0$, and $\al$ has the same predecessor $\gamma$. The key statement here is that the smallest such $\gamma > \beta$ is unique because we are on a simple tree $T$.  Now $\bV^\mu(\gamma) \le \bV^\mu (\al)\le 1$ as $\al\in \supp\mu$ and  potential $\bV$ of any positive measure on $T$ (and on  $T^2$) is a decreasing function always. 

\medskip

But  $\bV^\mu(\beta)=\bV^\mu(\gamma)$ because  $\bI^*(\tau) =0$ for all $\tau: \beta\le \tau <\gamma$ by the definition of $\gamma$ as the smallest interval containing interval $\beta$ for which $\mu(\gamma)>0$.

So we proved  that $\bV^\mu \le 1$ on $\supp\mu$ implies $\bV^\mu\le 1$ everywhere on $\pd T$. Then by monotonicity of potentials it is $\le 1$ everywhere on $T$.

\bigskip

This claim is blatantly false on $T^2$. The problems is that there can be a huge family $\Gamma$  of $\gamma>\beta$ such that $\mu(\gamma)>0$ and for any pair $\gamma_1, \gamma_2\in \Gamma$ none is smaller than the other. The reasoning above fails, and moreover there are plenty of simple examples of $\mu$ on $\pd T^2$ such that 
$$
\bV^\mu\le 1\,\,\text{on}\,\, \supp\mu,\,\,\text{but} \,\, \sup_{T^2}\bV^\mu\ge C,
$$
where $C$ is as large as one wishes (if $N$ is chosen large enough).

\bigskip

This phenomena is called {\it the lack of maximum principle}, and it reveals itself prominently in the following effect.

Let $\cT$ denote either $T$ or $T^2$.
Let us fix $\delta>0$ (not necessarily small but can be small) and consider
$$
E_\delta:=\{\al \in \cT: \bV^\mu(\al)\le \delta\}\,.
$$

Let 
$$
\bV_\delta^\mu(\al) = \bI(\one_{E_\delta} \bI^*\mu)(\al)\,.
$$ 
The expression (integration in the second equality is with respect to counting measure on $\cT$)
$$
\cE[\mu] =\int_{\cT} \bV^\mu\, d\mu= \int_{\cT} (\bI^*(\mu))^2
$$
is called the {\it energy of} $\mu$. The expression
$$
\cE_\delta[\mu] =\int_{\cT} \bV_\delta^\mu\, d\mu= \int_{E_\delta} (\bI^*(\mu))^2
$$
is called the {\it partial energy of} $\mu$. 

\bigskip

It is trivial that if $\cT=T$ then
\begin{equation}
\label{de}
\bV_\delta^\mu\le \delta
\end{equation}
uniformly. The reasoning is exactly the same as above for maximum principle. The consequence is the following partial energy estimate:
\begin{equation}
\label{en}
\cE_\delta[\mu] \le \delta\, |\mu|\,.
\end{equation}

\bigskip

But \eqref{de} can be easily false if $\cT=T^2$. We will show below the example that even \eqref{en} can be false.
All estimates in papers \cite{AMPS}, \cite{AMPVZ-K}, \cite{AHMV} are based on a weaker version of \eqref{en}, which  is true and which we call {\it the surrogate maximum principle}:
\begin{equation}
\label{smp0}
\cE_\delta[\mu] \le C \delta^\tau\, \cE[\mu]^{1-\tau}|\mu|^\tau\,.
\end{equation}
Here $C$ is universal.

For $\cT=T^3$ we can prove that with $\tau=1/3$, for $\cT=T^2$ we could originally prove it for $\tau =1/2$ and lately for any $\tau<1$. For $\tau=1$ it is false on $T^2$, see below. For $\cT=T^4$ we cannot prove \eqref{smp0}  at all, even for a very small $\tau$.

\section{Statement of the problem}
\label{state}

Suppose we are living on a rooted directed graph. For example on a dyadic tree $T$, or on $T^2=T\times T$. The latter can be viewed as a graph of all dyadic rectangles inside a unit square (the root). Let $\bI$ be operator of summation ``up the graph''. It has a formally adjoint operator $\bI^*$ of summation ``down the graph''. We use the same notation of this operator for the rooted dyadic tree $T$ and for  rooted bi-tree $T^2$. It is convenient to think that our graphs are finite, but very deep. The estimates and the constants must not depend on the depth. 
On dyadic tree $T$ we have the following key ``majorization theorem with small energy'':
\begin{theorem}
\label{d1}
Let $f, g: T\to \bR_+$, and  1) $g$ is superadditive, 2)  $\supp f \subset \{\bI g \le \delta\}$.  Let $\la \ge 10\delta$. Then there exists
$\vf: T\to \bR_+$  such that
\begin{enumerate}
\item $\bI \vf \ge \bI f$ on $\{2\la \le \bI g \le 4\la\}$;
\item $\int_T \vf^2 \le C\frac{\delta^2}{\la^2} \int_T f^2$.
\end{enumerate}
\end{theorem}

For a while we tried to prove the similar statement for $T^2$. Namely, we conjectured 
\begin{conj}
\label{d2}
Let $f, g: T^2\to \bR_+$, and  1) $g$ is superadditive in each variable, 2)  $\supp f \subset \{\bI g \le \delta\}$.  Let $\la \ge 10\delta$. Then there exists
$\vf: T^2\to \bR_+$  such that
\begin{enumerate}
\item $\bI \vf \ge \bI f$ on $\{2\la \le \bI g \le 4\la\}$;
\item $\int_{T^2} \vf^2 \le C\frac{\delta}{\la} \int_{T^2} f^2$.
\end{enumerate}
\end{conj}

For some very special cases, e. g. for $f=g$, this has been proved, and turned out to be a key result in describing the embedding  measures for the Dirichlet spaces in tri--disc into $L^2(\bD^3, d\rho)$. See
\cite{MPV}, \cite{AHMV}.

\section{Counterexample to small energy majorization on bi-tree}
\label{sem-cex}

Now we will show that this is not true in general. 

Moreover, below $f, g$ have special form, namely
$$
f=\bI^*\mu, \,g=\bI^*\nu,
$$
 with certain positive measures on $T^2$.  And measure $\mu$ is trivial, it is a delta measure of  mass $1$ at the root $o$ of $T^2$. In particular, $f(o)=1, f(v)=0, \forall  v\neq o$. Also $\bI f \equiv 1$ on $T^2$.

\bigskip

The choice of $\nu$ is more sophisticated. Choose large $\log n= 2^s$ and denote $2^M:= \frac{n}{\log n}$.

In the unit square $Q^0$ consider dyadic sub-squares $Q_1, \dots, Q_{2^M}$, which are South-West to North-East diagonal squares of size $2^{-M}\times 2^{-M}$.

In each $Q_j$ choose $\om_j$, the South-West corner dyadic square of size $2^{-n} \cdot 2^{-M}$.

Measure $\nu$ is the sum of delta measures at $\om_j, j=1, \dots, \frac{n}{\log n}$, each of muss $\frac1{n^2}$. Obviously 
$$
g(o)= \bV^* \nu(o) = \|\nu\|= \frac1{n^2}\cdot \frac{n}{\log n}= \frac1{n\log n}=:\delta.
$$
So we chose $\delta$ and $f, g$ satisfy $\supp f =\{o\}\subset \{ \bI g \le \delta\}$. Also $g$ is sub-additive in both variables on $T^2$: it is true for any function of the form $\bI^*\nu$.

\medskip

Now what is $\la$, and what is the set $\{2\la \le \bI g \le 4\la\}$?

\medskip

Consider (by symmetry this will be enough) $Q_1$ and $\om_1$ and consider the family $\cF_1$ of dyadic rectangles containing $\om_1$ and contained in $Q_1$ of the following sort:
$$
[0, 2^{-n} 2^{-M}]\times [0, 2^{-M}], [0, 2^{-n/2} 2^{-M}]\times [0, 2^{-2}2^{-M}], \dots, [0, 2^{-n/2^k} 2^{-M}]\times [0, 2^{-2^k}2^{-M}],
$$
there are approximately $\log n$ of them, and they are called $q_{10}, q_{11},\dots, q_{1k}$, $k\asymp \log n$. We do the same for each $\om_j, Q_j$ and we get $q_{j0}, q_{j1},\dots, q_{jk}$.

\begin{lemma}
\label{g}
$\bI g(q_{ji}) \asymp \frac1{n}\quad \forall j, i$.
\end{lemma}

It is proved in \cite{MPV}.

\bigskip

Let 
\begin{equation}
\label{F}
F:=\cup_{ik} q_{ik}\,.
\end{equation}

\bigskip

So we choose $\la=\frac{c}n$ with appropriate $c$.  Then
$$
F\subset \{2\la \le \bI g \le 4\la\}\,.
$$

As it was said that always $\bI f\ge 1$, so if $\vf$ as in Conjecture \ref{d2} would exist, we would have $\bI \vf \ge 1$ on $F$ and  (by the second claim of Conjecture \ref{d2})
$$
\int_{T^2} \vf^2 \le \frac{C}{\log n} \int_{T^2} f^2 = \frac{C}{\log n}\,.
$$
By the definition of capacity this would mean that
$$
\text{cap} (F) \le \frac{C}{\log n}\,.
$$

In the next subsection \ref{capFss} we show that $\text{cap} (F) \asymp 1$. Hence, conjecture is false.

\begin{remark}
\label{anyrate}
Moreover, this example shows that even a weaker estimate
$$
\int_{T^2} \vf^2 \le C\Big(\frac{\delta}{\la}\Big)^\eps \int_{T^2} f^2
$$
 is unattainable. Even more so, for any estimate of this type
$\int_{T^2} \vf^2 \le Ch\Big(\frac{\delta}{\la}\Big) \int_{T^2} f^2$, where $\lim_{t\to 0}h(t)=0$, we can use the construction above as a counterexample.
\end{remark}

\subsection{Capacity of $F$ is equivalent to $1$}
\label{capFss}

Let $\rho$ on $F$ be a capacitary measure of $F$, and $\rho_{jk}$ be its mass on $q_{jk}$.
By symmetry $\rho_{jk}$ does not depend on $j=1, \dots, \frac{n}{\log n}$.

\medskip 

The proof of the fact that 
$$
\rho_k:= \rho_{jk},\,\, j=1, \dots, \frac{n}{\log n}
$$
have the average $\ge \frac{c_0}{n}$, that is that
\begin{equation}
\label{E}
\bE \rho = \frac{\sum_{k=0}^{\log n} \rho_k}{\log n} \ge \frac{c_0}{n}
\end{equation}
follows below.

\bigskip  

In its turn it gives the required
\begin{equation}
\label{capF}
\text{cap} F \asymp 1\,.
\end{equation}

\medskip

Let us  first derive \eqref{capF} from \eqref{E}.
Measure $\mu$ that charges $\rho_k$ on each $q_{jk}, j=1,\dots, \frac{n}{\log n}; k= \frac{\log n}{2},\dots, \frac{3\log n}{4}$ is equilibrium so it gives $\bV^\mu \equiv 1$ on each $q_{jk}$.
Then \eqref{capF} follows like this: $\text{cap} F= \|\mu\|\ge \frac{n}{\log n} \sum_{k=0}^{\log n} \rho_k =  n \bE \rho$. Hence $\text{cap} F= \|\mu\|\ge  c_0$ if \eqref{E} is proved.

\medskip

Now let us prove \eqref{E}.
Everything is symmetric in $j$, so let $j=1$ and let us fix $k$ in $[\frac{\log n}{2}, \frac{3\log n}{4}]$. We know that
$$
1\le \bV^\mu \quad \text{on} \,\, q_{1k},
$$
and now let us estimate this potential from above.  For that we split $\bV^\mu$ to
$\bV_1$, this is the contribution of rectangles  containing $Q_1$, to
$\bV_2$, the contribution of rectangles containing $q_{1k}$ and contained in $Q_1$, and
$\bV_3$, the contribution of rectangles containing $q_{1k}$ that strictly intersect $Q_1$ and that are ``vertical'', meaning that there vertical side contains vertical side of $Q_1$.
(There is $\bV_4$ totally symmetric to $\bV_3$.)

Two of those are easy, $\bV_1$ ``almost'' consists of ``diagonal squares containing $Q_1$. Not quite, but other rectangles are also easy to take care. 
Denote
$$
r=\|\mu\|, \quad M=\log\frac{n}{\log n}\,.
$$
Then we write diagonal part first and then the rest:
$$
\bV_1= r+\frac{r}{2} +\frac{r}{4} + \dots \frac{r}{2^M}  +\frac{r}{2}  +\frac{r}{2} + 2\frac{r}{4} + 2\frac{r}{4} +
\dots   k\frac{r}{2^k}  + 2\frac{r}{2^k}  +\dots = C_1 r
$$

To estimate $\bV_2$ notice that there are at most $c n$ rectangles containing $q_{1k}$ and contained in $Q_1$ that do not contain any other $q$, there are $\frac{c n}{2}$ of rectangles contain $q_{1k}$ and one of its sibling (and lie in $Q_1$),  there are $\frac{c n}{4}$ of rectangles contain $q_{1k}$ and two of its sibling (and lie in $Q_1$), et cetera.

Hence,
$$
\bV_2 \le Cn\rho_k + \frac{Cn}{2} \rho_{k\pm 1} +  \frac{Cn}{4} \rho_{k\pm 2}+\dots
$$

Now consider $\bV_3$. The horizontal size of $q_{1k}$  is $2^{-M}\cdot 2^{-n2^{-k}}$. Its vertical size
is $2^{-M}\cdot 2^{-2^k}$. So the rectangles of the third type that do not contain the siblings: their number is
at most (we are using that $k\ge \frac12 \log n$)
$$
n2^{-k} ( 2^k+M) \le n + \sqrt{n}\log n\,.
$$
Those that contain $q_{1k}$ and one sibling, there number is at most
$$ 
n2^{-k} ( 2^{k-1}+M) \le \frac{n}2 + \sqrt{n}\log n\,.
$$
We continue, and get that
$$
\bV_3 \le n \rho_k + \frac{n}2  \rho_{k\pm1} + \frac{n}4 \rho_{k\pm 2}+\dots + \sqrt{n}\log n (\sum \rho_s)\,.
$$
Add all $\bV_i$:
$$
1 \le \bV_1+\bV_2 +\bV_3 +\bV_4 \le  C_1 r+ C n \rho_k+ \frac{Cn}2  \rho_{k\pm1} + \frac{Cn}4 \rho_{k\pm 2}+\dots + \sqrt{n}\log n (\sum \rho_s)\,.
$$
Now average over $k$. Notice that 
$$
r=\|\mu\|= \frac{n}{\log n} \sum\rho_s  =\frac12  n \bE \rho
$$
Hence,
$$
1 \le C' n \bE\rho + Cn \bE\rho + \frac{Cn}2 \bE\rho+  \frac{Cn}4 \bE\rho + \dots +  \frac12\sqrt{n}\log^2 n \,\bE\rho\,.
$$
Therefore, $\bE\rho\ge \frac{c_0}{n}$ and \eqref{E} is proved.

\section{The lack of $\int_{T^2} \bV^\nu_\eps \, d\nu \le C\eps |\nu|$ estimate}
\label{firstPower}

Let us recall the notations of \cite{MPV, AHMV, MPVZ}.   When we write $\int_{T^2}\dots$ without the indication of the measure of integrations, we mean the counting measure on $T^2$.
$$
\bV^\nu := \bI (\bI^*\nu)\,.
$$
$$
E_\eps:= \{(\al, \beta)\in T^2: \bV^\nu(\al, \beta) \le \eps\}\,.
$$
$$
\bV^\nu_\eps := \bI (\one_{E_\eps}\bI^*\nu)\,.
$$
$$
\cE[\nu]:= \int_{T^2} \bV^\nu\, d\nu= \int_{T^2} (\bI\nu)^2\,.
$$
This is energy of measure $\nu$.
$$
\cE_\eps[\nu]:= \int_{T^2} \bV^\nu_\eps\, d\nu= \int_{E_\eps} (\bI\nu)^2\,.
$$
This is a partial energy concentrated  where potential already became small.

\bigskip,

If we were on $T$ instead of $T^2$ we would have a trivial uniform estimate
\begin{equation}
\label{triv}
\bV^\nu_\eps \le \eps \Rightarrow  \cE_\eps[\nu]:= \int_{T} \bV^\nu_\eps\, d\nu \le C\eps |\nu|\,.
\end{equation}
And here $C=1$ of course.

The left hand side of \eqref{triv}, and what is more interesting, the right hand side of \eqref{triv} generally fails to hold on $T^2$ with any finite $C$. This section is devoted to a corresponding example.

\bigskip

On the other hand, \cite{MPVZ} proves that something like \eqref{triv} (but weaker) holds on $T^2$. In fact, we have  

\begin{theorem}[Surrogate Maximum Principle]
\label{smp}
If $\cE[\nu] \ge 2 \eps |\nu|$ then
$$
\cE_\eps[\nu] \le \eps e^{c_0\sqrt{\log\frac{\cE[\nu]}{\eps |\nu|}}}|\nu|\,.
$$
\end{theorem}

From Theorem \ref{smp} it is easy to deduce the following more transparent estimates:
\begin{equation}
\label{sqlog}
|\nu|\le \cE[\nu]\Rightarrow \cE_\eps[\nu] \le  C\eps e^{c\sqrt{\log\frac1{\eps}}} \,\cE[\nu]\,.
\end{equation}

\begin{equation}
\label{tau}
\cE_\eps[\nu]  \le C_\tau \eps^{1-\tau} \cE[\nu]^\tau |\nu|^{1-\tau}, \quad \forall \tau\in (0,1)\,.
\end{equation}

\bigskip

Let us build an example of sequences of pairs $\nu, \eps$  that shows that $\tau=0$ cannot be chosen in \eqref{tau}.

\bigskip

In other words let us build a counterexample to this ``universal'' estimate on $T^2$:
\begin{equation}
\label{fd}
``\int_{T^2} \bV_\eps^\nu\, d\nu \le C \eps |\nu|"\,.
\end{equation}

\begin{remark}
\label{FF}
One remark is in order. Notice that change of variables $\eps\to t\eps, \nu\to t\nu$ gets both the left hand side and the right hand side of \eqref{fd} multiplied by the same $t^2$. Thus we can normalize measure and always think that $|\nu|=1$.
Inequality above becomes $\int_{T^2} \bV_\eps^\nu\, d\nu \le C \eps$ for probability measures $\nu$. We show in this Section that it is false. But, in fact, in Section \ref{open} we show that {\it absolutely any} inequality
$$
\int_{T^2} \bV_x^\nu\, d\nu \le C F(x)
$$
is false regardless of function $F$. Notice that on $T$ function $F(x)=x$ makes the above inequality valid. So the counterexample in Section \ref{open} supersedes the counterexample we are going to explain in the current section. But in fact, a simple inspection shows that  both counterexamples are based on the same idea.
\end{remark}

To disprove \eqref{fd} we go back to the previous Sections and given $n$ (large power of $2$) we consider $\nu$ as in the previous Sections, and
$$
g=\bI^* \nu, \, \eps =\frac{c}{n},\, G:= g \cdot \one_{\bI g \le \frac{c}{n}}\,.
$$

Consider the set of vertices (rectangles) $F$ introduced in \eqref{F}. Lemma \ref{g}  claims that $\bI g\asymp \frac1n$ on $F$.
But then (with the right choice of $c$) $\one_{\bI g\le \frac{c}{n}} \bI g \asymp \frac1n$ on $F$. But it is easy to see that
$$
\bI g \cdot \one_{\bI g \le \la} \le \bI (g \cdot \one_{\bI g \le \la} )\,.
$$
Thus
\begin{equation}
\label{G}
\bI (G) = \bI (g\cdot \one_{\bI g \le \frac{c}{n}}) \ge \frac{c_0}{n}, \quad\text{on}\,\, F\,
\end{equation}

Suppose now that we want to bring to contradiction:
\begin{equation}
\label{contr}
\int_{T^2} G^2=\int \bV^\nu_{\frac{c}{n}} \, d \nu\le \frac{C}{n} |\nu|= \frac{C}{n} \cdot\frac1{n\log n}=\frac{C}{n^2\log n}\,.
\end{equation}

The definition of capacity and \eqref{G}, \eqref{contr} combined give us
$$
\text{cap}(F) \lesssim \frac1{\log n}\,.
$$

We come to a contradiction, because we proved in the previous Section that $\text{cap} (F) \asymp 1$. Contradictions shows that the only inequality in \eqref{contr} fails.

\section{The shape of the  graph of function $x\to \text{cap}( \bV^\nu \ge x)$}
\label{graph}

Below trees and bi-trees are finite, but unboundedly deep.
Let $E$ be a subset of $T$ or $T^2$ and $\nu$ be a capacitary measure for $E$, 
$$
\text{cap}(E) =|\nu|, \, \bV^\nu=1 \,\, \text{on} \,\, \supp\nu, \,\, g:= \bI^* \nu=\{f: \int f^2\to \min\,\, \text{for}\,\, \bI f \ge 1\,\, \text{on}\,\, E\}\,.
$$

\bigskip

First consider the case of $T$.  Let $x\in [|\nu|, 1]$ and we study the set 
$$
D_x:=\{\al\in T: \bV^\nu(\al)\ge x\}\,.
$$
We want to understand a bit the shape of the graph of
$$
C(x) :=\text{cap}(D_x)\,.
$$
We start with $x=|\nu|=\text{cap}(E)$. Notice that $o$, the root of $T$, is obviously such that $\bV^\nu(o) =|\nu|$, so $0\in D_{|\nu|}$. But $\text{cap}(0) = \text{cap}(T)=1$. Thus
$$
C(|\nu|) =1\,.
$$
Now consider $x=1$. On $E$ we have $\bV^\nu=1$ and maximum principle (we are on $T$, so it exists) says that $E=\{\al: \bV^\nu \ge 1\}$. Therefore,
$$
C(1) = \text{cap}(E) =|\nu|\,.
$$

Now let $|\nu|<x<1$.
We  know (again this is maximum principle) that
\begin{equation}
\label{MP}
\int_T \one_{\bI g \le x} \cdot g^2 = \int_T \bV^\nu_x \, d\nu \le x |\nu|\,.
\end{equation}

Notice that if $\bI g(\al)\le x$ and $\bI g (\text{son}\, \al)>x$ then $\bI g(\al)\ge x/2$ just because $g=\bI^*\nu$ is monotonically increasing on $T$. But this means that
\begin{equation}
\label{cut}
\bI ( \one_{\bI g \le x} \cdot g)  \ge x/2, \quad \text{on} \,\, D_x=\{ \bI g = \bV^\nu \ge x\}\,.
\end{equation}

The definition of capacity and relationships \eqref{MP}, \eqref{cut} show the following:
\begin{theorem}
\label{capT}
On a simple tree $T$ the capacity of the level set $D_x=\{\al\in T: \bV^\nu(\al)\ge x\}$ for any capacitary measure $\nu$ of  a set $E$ satisfies the following inequality
$$
C(x)=\text{cap}(\{\al\in T: \bV^\nu(\al)\ge x\}) \le \frac{4\text{cap}(E)}{x}=  \frac{4|\nu|}{x},\,\, \text{cap}(E)\le  x\le 1\,.
$$
\end{theorem}

\bigskip

This is absolutely not the case for $T^2$. The capacity of level set of capacitary potentials on $T^2$ behave in a much stranger and wild way. We saw it in Section \ref{capF}. In fact, our measure $\nu$ in the previous Sections is a capacitary measure, 
$$
|\nu| =\frac1{n\log n}\,.
$$
We put 
$$
x=\frac{c}n\,.
$$
 We saw in Section \ref{capF} that if the absolute constant $c$ is chosen correctly, then
\begin{equation}
\label{capT2}
\text{cap}((\al, \beta) \in T^2: \bV^\nu(\al, \beta) \ge \frac{c}n) \asymp 1>> \frac{|\nu|}{x}\,.
\end{equation}
This means that Theorem \ref{capT} is false for $T^2$ because if it were true, that we would have $\text{cap}((\al, \beta) \in T^2: \bV^\nu(\al, \beta) \ge \frac{c}n)\lesssim \frac1{\log n}$.

\bigskip

\subsection {The reason for the effect \eqref{capT2}}
\label{reason}

On $T^2$ we do not have \eqref{triv}, which is \eqref{MP} above. Instead we have \eqref{tau} that makes the estimate of capacity much faster blowing up than in Theorem \ref{capT}. In fact, \eqref{tau} claims
$$
\text{cap}(\{\bV^\nu\ge x\}) \le \frac{C_\tau\text{cap}(E)}{x^{1+\tau}}\,.
$$
and we saw that $\tau$ is indispensable. Of course the capacity of any subset of $T^2$ is bounded by $1$, so we have
$$
\text{cap}(\{\bV^\nu\ge x\}) \le \max\Big( 1, \frac{C_\tau\text{cap}(E)}{x^{1+\tau}}\Big)\,.
$$
This explains a flat piece of graph $C(x)\asymp 1$, when $x$ is between $\frac1{n\log n}$ and $\frac1n$.
In fact, looking at \eqref{sqlog} we may think that the flat piece of graph can be much wider.

\section{One more counterexample}
\label{open}
Here is the question asked by Fedor Nazarov. He also hinted us a possible  construction of a counterexample. 

\noindent{\bf Question.} Consider normalized measures on the unit square, $|\mu|=1$. Let $x>>1$. Is it always possible to have the estimate
$$
\int_{T^2} \bV^\mu_x \, d\mu=\int_{\bV^\mu\le x} (\bI^*\mu)^2 \le F(x)\,?
$$

\medskip

The meaning of this question is that we always (see Theorem \ref{smp}, and   \eqref{sqlog}, \eqref{tau}) have some trace of total energy
in the right hand side of our estimates of partial energy. What if total energy is huge or   ``infinite''? Maybe one does not need this total energy contribution into the right hand side as its presence in Theorem \ref{smp}, and in \eqref{sqlog}, \eqref{tau})? Maybe the partial energy is always bounded by a function of its ``cut-off'' parameter $x$ for all normalized measures?

We will show that no estimate as above exists (but on $T$ it does exist with the simplest $F(x)=x$).

\bigskip

Let us fix two large integers $n, M$ ($n$ will be much bigger than $M$) and consider a small modification of the construction
of previous Sections. Namely,
Consider $2^M$ dyadic squares located on the diagonal of $Q^0:=[0,1]^2$, each of size $2^{-M}\times 2^{-M}$. We call them $Q_1, \dots, Q_{2^M}$.

In South West corner of each $Q_j$ choose a dyadic square of size $2^{-n-M}\times 2^{-n-M}$. Call them $\om_1, \dots, \om_{2^M}$. We charge each $\om_j$ with mass $2^{-M}$, this forms measure $\mu$ of mass $1$. Now consider  $\om_1$ and
the family $\cF_1$ of dyadic rectangles containing $\om_1$ and contained in $Q_1$ of the following sort:
$$
[0, 2^{-n} 2^{-M}]\times [0, 2^{-M}], [0, 2^{-n/2} 2^{-M}]\times [0, 2^{-2}2^{-M}], \dots, [0, 2^{-n/2^k} 2^{-M}]\times [0, 2^{-2^k}2^{-M}],
$$
there are approximately $\log n$ of them, and they are called $q_{10}, q_{11},\dots, q_{1k}$, $k\approx \log n$.
We repeat this for other $\om_j$, having now $q_{j0}, q_{j1},\dots, q_{jk}$, $j=1, \dots, 2^M$, $k\approx \log n$, and $q_ji$ containing $\om_j$ and contained in $Q_j$.

Now put
$$
x= n 2^{-M}
$$
and choose $n$ to have $x>>1$.

\medskip

We claim that 
\begin{equation}
\label{Vx}
\bV(q_{ji}) \le Cx\quad \forall j, i\,.
\end{equation}

We know from \cite{MPV,MPVZ} that given $j, i$ there are approximately $n$ dyadic rectangles containing $q_{ji}$ and contained on $Q_j$.  Each gives contribution $2^{-M}$ into $\bV(q_{ji})$. So if we would count only them in $\bV(q_{ji})$ then we get the total of $\approx n2^{-M}$, and \eqref{Vx} would follow. Let us call this contribution of $\approx n2^{-M}$ obtained as above {\bf the main contribution}. Let us justify that it is the main contribution now. 

But there are much more dyadic rectangles containing $q_{ji}$ and contained in $Q^0$. Let us bookkeep their contributions to 
$\bV_x(q_{ji})$. We hope that those are not too big in order to prove \eqref{Vx}.

Notice that if  \eqref{Vx} is proved we have many rectangles $R$ with  $\bV(R) \le Cx$; so many that we can hope to prove that 
\begin{equation}
\label{large}
\sum_{R: \bV^\mu(R) \le Cx} \mu(R)^2 \ge F(x)\,.
\end{equation}

\medskip

So we fix, say,  $q_{0i} =[0, 2^{-n/2^i} 2^{-M}]\times [0, 2^{-2^i}2^{-M}]$, and we can see that apart of $\approx n$ rectangles between $q_{0i}$ and $Q_0$, there are also

a) tall rectangles $[0, 2^{-n/2^{i'}} 2^{-M}]\times [0, 2^m 2^{-M}]$, $i\le i' \le \log n$, $1\le m\le M$, containing $q_{0i}$;

b) long rectangles $[0, 2^{m} 2^{-M}]\times [0, 2^{-2^{j'}}2^{-M}]$, $0\le j' \le j$, $1\le m\le M$, containing $q_{0i}$;

c) $m$-large rectangles, containing $q_{0i}$: these are rectangles containing dyadic square $Q_0^{(m)}$ with side $2^m 2^{-M}$
that contains $Q_0$, but not containing $Q_0^{(m+1)}$, $m=2, \dots, M$.

The contribution of tall rectangles into $\bV(q_{0i})$ is bounded by $M2^{-M} \log n$, and the same is for the contribution of the 
long rectangles. 

The contribution of $M$-large rectangles is $1$--in fact, there is only one such rectangle, namely our initial unit square $Q^0$. 
The contribution of $M-1$-large rectangles is $\frac12 \cdot (1+1+1)$. In fact we would have $3$ rectangles in the family of $M-1$-large rectangles: $Q_0^{M-1}$ square itself and its two predecessors, one long, one tall.

The contribution of $M-2$-large rectangles is $\frac14 \cdot (1+2+2)$. In fact we would have $5$ rectangles in the family of $M-2$-large rectangles: $Q_0^{M-2}$ square itself and its four predecessors, two long, two tall.

The contribution of $M-3$-large rectangles is $\frac18 \cdot (1+3+3)$. In fact we would have $5$ rectangles in the family of $M-2$-large rectangles: $Q_0^{M-3}$ square itself and its six predecessors, three long, three tall.

Et cetera.  The total contribution of all $m$-large rectangles containing $q_{0i}$ is at most
$$
\sum_{m=1}^M \frac1{2^{m}} (2m+1)\le C_1\,.
$$
This is definitely smaller that $x=n2^{-M}$ than then main contribution and can be just absorbed into the main contribution $\approx n2^{-M}=x>>1$.

The contributions from long and tall rectangles listed in a) and b) above is at most $2\log n M2^{-M} \lesssim n^\tau 2^{-\tau M} \lesssim (n 2^{-M})^\tau =x^\tau<< x$ for any $\tau>0$ for example  for $\tau=\frac12$. Hence, the contribution from  long and tall rectangles listed in a) and b) above is also much smaller than the main contribution of order $x$ and can be absorbed into the main contribution. 

We finally proved \eqref{Vx}. Now let us estimate $\sum_{R: \bV^\mu(R) \le Cx} \mu(R)^2 $ from below. From \cite{MPV, MPVZ} we know that for each $q_{ji}$ there is a family of dyadic rectangles $\cF_{ji}$ such that 1) every $R\in \cF_{ji}$ contains $q_{ji}$ and is contained in $Q_j$, $j=1, \dots, 2^M$,  2) the cardinality of $\cF_{ji}$ is at least $c\,n$, $c>0$, 3) families $\cF_{ji}$ are disjoint, $j=1, \dots, 2^M$, $i\le C\log n$. Each rectangle $R$ of $\cup_j\cup_i \cF_{ji}$ has the property that
$$
\bV^\mu(R) \le Cx\,.
$$
We proved this in \eqref{Vx}. So each of those $R$ gives a contribution into the sum $\sum_{R: \bV^\mu(R) \le Cx} \mu(R)^2 $, and this contribution is $2^{-2M}$. Therefore,
$$
\sum_{R: \bV^\mu(R) \le Cx} \mu(R)^2 \ge 2^{-2M}\cdot \sharp j \cdot \sharp i \cdot \sharp (\cF_{ji}) \ge c\, 2^{-2M} 2^M \log n \cdot n =c\,2^{-M} n \cdot  \log n=
$$
$$
c\, x \cdot (\log x +M)\,.
$$
Now, given $x>>1$, we can freely choose $M$, e.g. $M=x, x^2, 2^x, F(x)\dots$ and then choose $n$ from $n2^M=x$ and do the construction above. So \eqref{large} is proved.


\end{document}